\documentclass[journal]{IEEEtran}
\ifCLASSINFOpdf
  \usepackage[pdftex]{graphicx}
\else
\fi
\usepackage{amsmath}
\usepackage{amsfonts}
\usepackage{symbols}
\def\figuref{Fig.~}

\def\tableref{TABLE~}

\begin{document}
%
\title{A discretisation method with the $H_{\rm div}$ inner product for
electric field integral equations}
%
%
%

\author{Kazuki~Niino,~\IEEEmembership{Member,~IEEE,}
        Sho~Akagi\IEEEmembership{}
        and~Naoshi~Nishimura\IEEEmembership{}
\thanks{K. Niino, S. Akagi and N. Nishimura are with the Kyoto University, Kyoto, Japan.}
}

\maketitle

\begin{abstract}
A discretisation method with the $H_{\rm div}$ inner product for the electric field integral equation~(EFIE) is proposed. 
The EFIE with the conventional Galerkin discretisation shows bad accuracy for problems with a small frequency, a problem known as the low-frequency breakdown.
The discretisation method proposed in this paper utilises the $H_{\rm div}$ scalar product with a scalar coefficient for the Galerkin discretisation
and overcomes the low-frequency problem with an appropriately chosen coefficient.
As regards the preconditioning,
we find that a naive use of the widely-used Calderon preconditioning is not efficient for reducing the computational time with the new discretisation.
We therefore propose a new preconditioning which can accelerate the computation successfully.
The efficiency of the proposed discretisation and preconditioning is verified through some numerical examples.
\end{abstract}


%
\IEEEpeerreviewmaketitle

\section{Introduction}
The boundary element method (BEM), which is also called the method of moment (MoM) in electromagnetic community, 
is one of well-known methods for solving electromagnetic problems.
Various formulations of boundary integral equations for EM applications have been proposed,
among which is the electric field integral equation (EFIE) \cite{chew1995waves} which 
is effective for scattering problems with perfect electric conductors (PECs).
It is known, however, that the EFIE suffers from bad accuracy when the frequencies are small (\cite{mautz1984, vecchi1999loop}).
This problem, called ``low-frequency breakdown'', is due to the ill-conditioning
of the coefficient matrix obtained by discretising the EFIE.
Indeed, some parts of discretised EFIE are lost when $k h\rightarrow 0$ where $k$ is the wave number and  $h$ is the average diameter of the mesh.
A widely used solution to the low-frequency breakdown is the loop-tree decomposition \cite{wu1995study, vecchi1999loop}.
This method divides a discretised integral equation into two sets of linear equations with the help of the quasi--Helmholtz decomposition and rescales these linear equations so that they do not vanish when $k h \rightarrow 0$.
Another solution to the low-frequency breakdown is the augmented integral equation \cite{qian2008augmented}.
This method solves the current continuity equation simultaneously with the standard EFIE
with the surface electric charge as additional unknowns.
Both methods can remedy the low-frequency breakdown, but the additional computational time introduced by the loop-star decomposition or the new set of equations and unkonwns is not ignorable.
We have found that the low-frequency breakdown can also be avoided as one uses
the $H_{\rm div}$ inner product for the Galerkin method instead of the $L^2$ inner product \cite{niino2013maxwelljp}.
This method reduces the conventional discretised integral equation to a weighted sum of itself and its surface divergence.
We have verified that this method can remedy the low-frequency breakdown as one chooses an appropriate constant for the $H_{\rm div}$ inner product in the Poggio-Miller-Chang-Harrington-Wu-Tsai (PMCHWT) formulation \cite{niino2013maxwelljp}.

The EFIE also has the problem of slow convergence
when it is solved with an iterative linear solver such as the generalised minimal residual method (GMRES) \cite{saad2003}.
This problem occurs since the electric field integral operator (EFIO) is ill-conditioned.
The convergence of the EFIE becomes worse as a finer mesh is used
since the condition number of the discretised EFIE is proportional to $1/h^2$\cite{andriulli2008multiplicative}.
Hence, an acceleration of the iteration method, typically Calderon's preconditioning, is indispensable with the EFIE.
The Calderon preconditioning was first proposed by Steinbach and Wendland for Laplace's equation \cite{steinbach1998construction}
and was applied to the EFIE by Christiansen and Nedelec \cite{christiansen2003preconditioner}.
A multiplicative Calderon preconditioning can be constructed \cite{andriulli2008multiplicative} with the help of the Rao-Wilton-Glisson (RWG) basis function \cite{rao1982electromagnetic} and the Buffa-Christiansen (BC) basis functions \cite{buffa2005dual}.
However, standard EFIEs with Calderon's preconditioning still suffer from the low-frequency breakdown.

In this paper, we propose a preconditioned EFIE discretised with the $H_{\rm div}$ inner product,
which solves both the low-frequency breakdown and ill-conditioning.
We found that a naive use of the Calderon preconditioning in the EFIE discretised with the $H_{\rm div}$ inner product cannot reduce the computational time efficiently.
We, therefore, introduce another preconditioning method which does decrease the computational time.

The additional computational time of the proposed method for solving the low-frequency breakdown is small, compared with conventional methods
such as the loop-star decomposition and the method of the augmented integral equation.
In fact, the proposed method requires
the calculation of the normal component of the MFIE in addition to the standard EFIE.
But, the additional computational time for calculating the MFIE is small with the fast multipole method (FMM)
since the most parts of the FMM computation are common to the EFIE and the MFIE.

This paper is organised as follows.
In section \ref{sec:formulation},
we formulate the electromagnetic wave scattering problems and the EFIE.
In section \ref{sec:disc_conv}, we introduce the conventional discretisation method
and the low-frequency breakdown.
Then, we propose a discretisation with the $H_{\rm div}$ scalar product and 
describe how this method solves the low-frequency breakdown in section \ref{sec:disc_proposed}.
We introduce an effective multiplicative preconditioning to the proposed method in section \ref{sec:calderon}.
After this, we show the effectiveness of the proposed method via some numerical examples in section \ref{sec:num_ex}
and make conclusion in section \ref{sec:conclusion}.

\section{Formulation}\label{sec:formulation}
\begin{figure}[htbp]
  \centering
      \includegraphics[width=4cm]{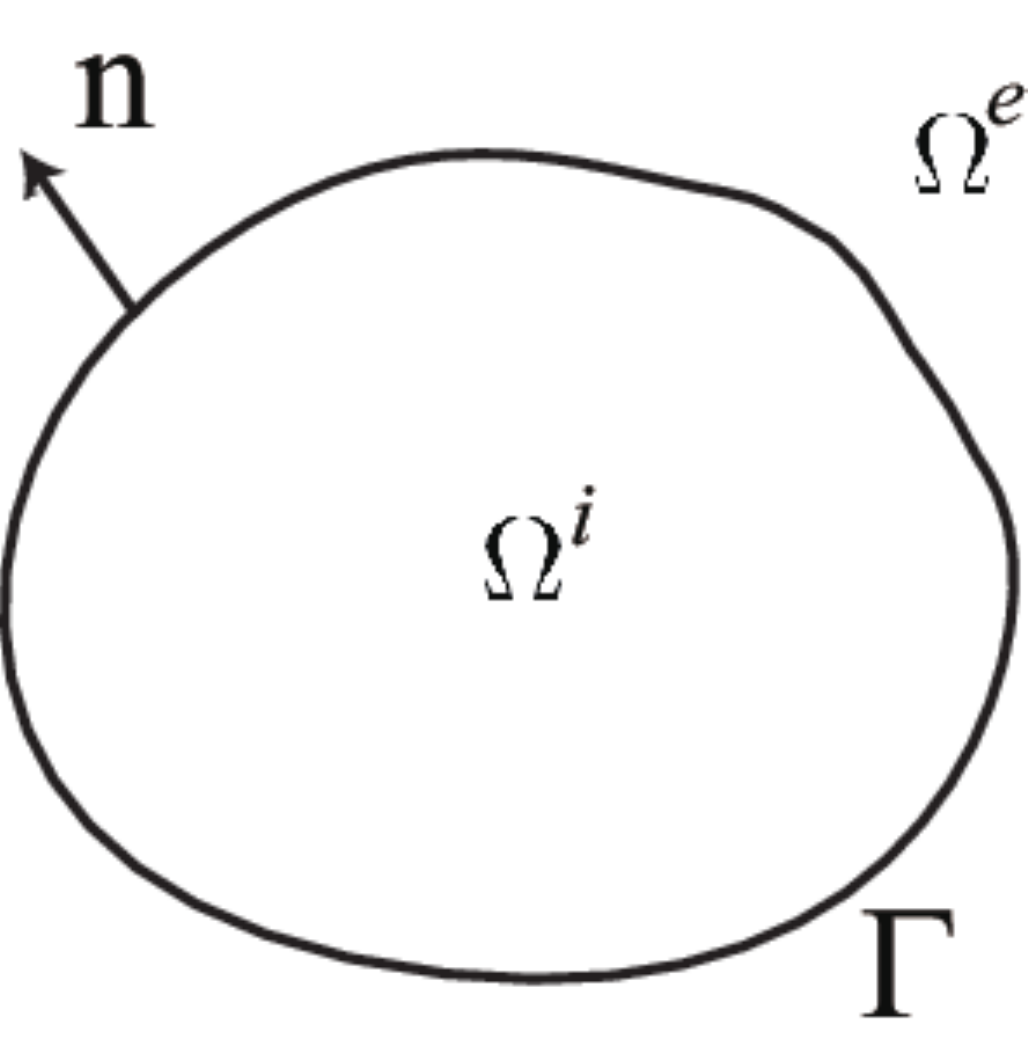}
    \caption{Electromagnetic scattering problem.}
    \label{efie_dom}
\end{figure}

We consider electromagnetic wave scattering problems with a single PEC as shown in \figuref\ref{efie_dom}.

The domain of the PEC is denoted by $\Omega^i$ and is enclosed by the smooth boundary $\Gamma$.
We find the solutions $\bE$ and $\bH$ satisfying the Maxwell equations
\begin{align*}
 \nabla\times\bE = \zi\omega\mu \bH,\quad  \nabla\times\bH = -\zi\omega\varepsilon \bE
\end{align*}
in $\Omega^e=\mathbb{R}^3\backslash\overline{\Omega^i}$, the boundary condition
\begin{align*}
 \bE^+\times\bn=0
\end{align*}
on $\Gamma$ and the radiation conditions for scattered waves $\bE^{\rm sca}$ and $\bH^{\rm sca}$ where $\bE$ and $\bH$ are unknown electric and magnetic fields,
$\omega$ is the frequency
with the time dependency of various quantities being $\ex^{-\zi\omega t}$,
$\varepsilon$ and $\mu$ are the permittivity and permeability in $\Omega^e$,
$\bn$ is the exterior unit normal to $\Gamma$,
$\bE^+$ is the limit value of $\bE$ from $\Omega^e$ to $\Gamma$,
$(\bE^{\rm sca}, \bH^{\rm sca})=(\bE-\bE^{\rm inc}, \bH-\bH^{\rm inc})$
and $\bE^{\rm inc}$ and $\bH^{\rm inc}$ are electric and magnetic incident waves, respectively.

The EFIE for this problem can be written as follows:
\begin{align}
\label{eq:efie}
 \zi\omega\mu Q \bj = \bE^{\rm inc}\times\bn,
\end{align}
where $\bj$ is the unknown electric current on $\Gamma$,
\begin{align}\label{eq:tau}
Q\bj =  \bn\times\int_\Gamma\left\{G(\bx-\by) 
   + \frac{1}{k^2}\nabla\nabla G(\bx-\by)\right\}\bj(\by)\md S_y,
\end{align}
$k=\omega\sqrt{\varepsilon\mu}$ and $G$ is Green's function of the Helmholtz equation:
\begin{align*}
 G(\bx-\by) = \frac{\ex^{\zi k|\bx-\by|}}{4\pi |\bx-\by|}
\end{align*}

\section{Conventional Galerkin Method and Low-Frequency Breakdown}\label{sec:disc_conv}
In this section, we describe the conventional Galerkin discretisation method
for (\ref{eq:efie})
and show that the solution of the linear equations obtained in this way
may have a large error due to the low-frequency breakdown.

\subsection{Galerkin Method with the $L^2$ Inner Product}
In the conventional Galerkin method, (\ref{eq:efie}) is tested with the basis functions $\bt_i$ in the following way:
\begin{align*}
&\biggl(\bn\times\bt_i,\\ &\bn\times\int_\Gamma\biggl\{\zi\omega\mu G(\bx-\by)\bj (\by) 
   + \frac{\zi}{\omega\varepsilon}\nabla\nabla G(\bx-\by)\bj (\by)\biggr\}\md S_y\biggr)_{L^2_T(\Gamma)}\\
   = &(\bn\times\bt_i, \bE^{\rm inc}\times\bn)_{L^2_T(\Gamma)}.
\end{align*}
where $(\cdot, \cdot)_{L^2_T(\Gamma)}$ is the $L^2$ inner product of tangent vectors on $\Gamma$:
\begin{align*}
 (\bu, \bv)_{L^2_T(\Gamma)} = \int_\Gamma \overline{\bu}\cdot \bv \md S.
\end{align*}
Expanding the unknown function $\bj$ with
\begin{align*}
 \bj = \sum_{i=1}^N j_i \bt_i,
\end{align*}
we obtain the following linear equation
\begin{align}\label{eq:efie_l2}
 A_{L^2}\bx = \bb_{L^2},
\end{align}
where $\bt_i$ is the RWG basis function and $N$ is the number of the RWG basis functions.
The elements of the matrix $A_{L^2}$ and the vector $\bb_{L^2}$ are defined by
\begin{align*}
 (A_{L^2})_{ij} &= \biggl(\bn\times\bt_i,\bn\times\int_\Gamma\biggl\{\zi\omega\mu G(\bx-\by)\bt_j (\by) \\
   &+ \frac{\zi}{\omega\varepsilon}\nabla\nabla G(\bx-\by)\bt_j (\by)\biggr\}\md S_y\biggr)_{L^2_T(\Gamma)},\\
 (\bb_{L^2})_{i} &= (\bn\times\bt_i, \bE^{\rm inc}\times\bn)_{L^2_T(\Gamma)},
\end{align*}
where $(X)_{ij}$ is the $(i,j)$ element of a matrix $X$ and $(Y)_i$ is the $i$th element of a vector $Y$.
\subsection{Low-Frequency Breakdown}\label{sec:low_freq_conv}
Equation (\ref{eq:efie_l2}) is known to suffer from the low-frequency breakdown.
We now show that the solution of (\ref{eq:efie_l2}) obtained with iteration methods
may have a large error.
See Zhao and Chew \cite{zhao2000integral} for a related discussion for small $k$.


We first use the loop-star basis functions \cite{vecchi1999loop} for expanding the unknown function $\bj$.
The loop functions $\bt_m^{\rm loop}$ are defined by
\begin{align}\label{eq:def_p}
 \bt_m^{\rm loop} = {\rm curl}_S\, p_m, \\\notag
 {\rm curl}_S\, p_m := \bn\times\nabla p_m
\end{align}
where $p_m$ is the piecewise linear function associated with the vertex $m$,
which is $1$ at the vertex $m$ and decreases linearly to $0$ at neighboring vertices.
The star functions $\bt_n^{\rm star}$ are defined by the linear combination
of RWG basis functions associated with the three edges of the $n$th element \cite{wu1995study}:
\begin{align*}
 \bt_n^{\rm star} = \sum_{i=1}^3 \frac{S_{n_i}}{l_{n_i}} \bt_{n_i}
\end{align*}
where $n_i$ is the index of three edges of the $n$th triangle,
$l_{n_i} $ is the length of the edge $n_i$ and
$S_{n_i}$ is either $1$ or $-1$ which is defined so that
the current $S_{n_i}\bt_{n_i}$ flows out from the $n$th element.
It is known that the span of loop-star basis functions is identical with that of RWG basis functions \cite{vecchi1999loop}:
\begin{align}\notag
 \bj &=  \sum_{i=1}^N j_i\bt_i\\\label{eq:lt_decom}
  &=\sum_{i=1}^{N_{\rm loop}^{\rm RWG}} j_i^{\rm loop}\bt_i^{\rm loop}
  + \sum_{i=1}^{N_{\rm star}^{\rm RWG}} j_i^{\rm star}\bt_i^{\rm star}.
\end{align}
where $N_{\rm loop}^{\rm RWG}$ and $N_{\rm star}^{\rm RWG}$ are the number of independent loop and star functions of the RWG basis functions, respectively.
Note that we use the loop-star decomposition only for studying numerical methods,
but never in computations in this paper.

Substituting (\ref{eq:lt_decom}) into (\ref{eq:efie_l2}),
we obtain
\begin{align}\label{eq:efie_l2_decomposed}
 \begin{pmatrix}
  Z^{LL}_{L^2} & Z^{LS}_{L^2}\\
  Z^{SL}_{L^2} & Z^{SS}_{L^2}
 \end{pmatrix}
\begin{pmatrix}
 \bj^{\rm loop} \\ \bj^{\rm star}
\end{pmatrix}=
 \begin{pmatrix}
  \bb^L_{L^2} \\ \bb^S_{L^2}
 \end{pmatrix}
\end{align}
where $Z^{LL}_{L^2}, Z^{LS}_{L^2}, Z^{SL}_{L^2}$ and $ Z^{SS}_{L^2}$
 are matrices defined by
\begin{align*}
 (Z^{LL}_{L^2})_{ij} &= (\bn\times\bt^{\rm loop}_i, \zi\omega\mu Q\bt^{\rm loop}_j)_{L^2_T(\Gamma)}\\
 (Z^{LS}_{L^2})_{ij} &= (\bn\times\bt^{\rm loop}_i, \zi\omega\mu Q\bt^{\rm star}_j)_{L^2_T(\Gamma)},\\
 (Z^{SL}_{L^2})_{ij} &= (\bn\times\bt^{\rm star}_i, \zi\omega\mu Q\bt^{\rm loop}_j)_{L^2_T(\Gamma)},\\
  (Z^{SS}_{L^2})_{ij} &= (\bn\times\bt^{\rm star}_i, \zi\omega\mu Q\bt^{\rm star}_j)_{L^2_T(\Gamma)},
\end{align*}
and $\bj^{\rm loop}, \bj^{\rm star}, \bb^L_{L^2}$ and $ \bb^S_{L^2}$ are vectors defined as follows:
\begin{align*}
 (\bb^L_{L^2})_i &= (\bn\times\bt^{\rm loop}_i, \bE^{\rm inc}\times\bn)_{L^2_T(\Gamma)},\\
 (\bb^S_{L^2})_i &= (\bn\times\bt^{\rm star}_i, \bE^{\rm inc}\times\bn)_{L^2_T(\Gamma)},\\
 \bj^{\rm loop}&=(j^{\rm loop}_1, j^{\rm loop}_2, \cdots , j^{\rm loop}_{N_{\rm loop}^{\rm RWG}})^T,\\
 \bj^{\rm star}&=(j^{\rm star}_1, j^{\rm star}_2, \cdots , j^{\rm star}_{N_{\rm star}^{\rm RWG}})^T.
\end{align*}

Now, we estimate the order of each element in (\ref{eq:efie_l2_decomposed}) with respect to $k$ and $h$ under the condition that $k h$ is sufficiently small, where $h$ is the largest diameter of the triangular mesh.
The orders of the matrices $Z^{LL}_{L^2}, Z^{LS}_{L^2}, Z^{SL}_{L^2}$ and $ Z^{SS}_{L^2}$
are those of elements having the maximum absolute values.
Hence the orders of the matrices $Z^{LL}_{L^2}$ and $Z^{SS}_{L^2}$ are
equal to those of diagonal elements of these matrices,
and the order of $Z^{LS}_{L^2}$ ($Z^{SL}_{L^2}$) is that of one of the $ij$ elements where $\bt^{\rm loop}_i$ ($\bt^{\rm star}_i$) and $\bt^{\rm star}_j$ ($\bt^{\rm loop}_j$) share their supports.
In the following evaluation, we assume that the basis functions $\bt^{\rm loop}_i$ and $\bt^{\rm star}_i$ are normalised such that ${\rm sup}\,|\bt^{\rm loop}_i|={\rm sup}\,|\bt^{\rm star}_i|=1$.

The matrix $Z^{SS}_{L^2}$ satisfies
\begin{align}\notag
 Z^{SS}_{L^2} &= \left(\bn\times\bt_i^{\rm star}, 
 \zi\omega\mu \bn\times \int_\Gamma G\bt_i^{\rm star}\md S_y\right)_{L^2_T(\Gamma)}\\\notag
 &+\left(\bn\times\bt_i^{\rm star}, 
 \zi\omega\mu \bn\times \int_\Gamma \frac{1}{k^2}\nabla\nabla G\bt_i^{\rm star}(\by)\md S_y\right)_{L^2_T(\Gamma)}\\\notag
 &= \left(\bt_i^{\rm star}, 
 \zi\omega\mu \int_\Gamma G\bt_i^{\rm star}\md S_y\right)_{L^2_T(\Gamma)}\\\label{eq:zss}
 &-\left(\nabla_S\cdot\bt_i^{\rm star},\frac{\zi}{\omega\varepsilon} \int_\Gamma G \nabla_S\cdot \bt_i^{\rm star}\right)_{L^2_T(\Gamma)}
\end{align}
where
\begin{align*}
 \nabla_S \cdot\bphi &:= -(\nabla\times(\bphi\times\bn))\cdot \bn.
\end{align*}
These two terms in (\ref{eq:zss}) are estimated as
\begin{align*}
 \left(\bt_i^{\rm star}, 
 \zi\omega\mu \int_\Gamma G\bt_i^{\rm star}\md S_y\right)_{L^2_T(\Gamma)}
 &\sim O(k h^3),\\
 -\left(\nabla_S\cdot\bt_i^{\rm star},\frac{\zi}{\omega\varepsilon} \int_\Gamma G \nabla_S\cdot \bt_i^{\rm star}\right)_{L^2_T(\Gamma)}
 &\sim O\left(\frac{h}{k}\right)
\end{align*}
since 
\begin{align*}
 \bt_i^{\rm star}\sim O(1), \quad \nabla_S\cdot \bt_i^{\rm star} \sim O\left(\frac{1}{h}\right),\quad \int_{\Gamma_i} G\md S_y \sim O(h)
\end{align*}
and the area of a triangle is $O(h^2)$ where $\Gamma_i$ is the support of the basis function $\bt_i$.
Hence we obtain
\begin{align*}
 Z^{SS}_{L^2}\sim O\left(\frac{h}{k}\right).
\end{align*}
since $|k h^3| \ll |h/k|$ if $|k h| \ll 1$.
In $Z^{LL}_{L^2}, Z^{LS}_{L^2}$ and $ Z^{SL}_{L^2}$, however, the hyper-singular term vanishes
and the first term in the RHS of (\ref{eq:tau}) is dominant.
We, therefore, obtain
\begin{align*}
 Z^{LL}_{L^2} &\sim \left(\bn\times\bt_i^{\rm loop},
 \zi\omega\mu\bn\times\int_\Gamma G\bt^{\rm loop}_i \md S_y\right)\\
 &\sim O(k h^3)
\end{align*}
The same calculation can be applied to $Z^{LS}$ and $Z^{SL}$ and,
finally, we obtain the following evaluations:
\begin{align*}
 Z^{LL}_{L^2}\sim O(k h^3),\quad Z^{LS}_{L^2}\sim O(k h^3),\\
 Z^{SL}_{L^2}\sim O(k h^3),\quad Z^{SS}_{L^2}\sim O\left(\frac{h}{k}\right).
\end{align*}
For the RHS of (\ref{eq:efie_l2_decomposed}), we have
\begin{align*}
 (\bb^L_{L^2})_i &= (\bn\times\bt^{\rm loop}_i, \bE^{\rm inc}\times\bn)_{L^2_T(\Gamma)}\\
 &= -(\bt^{\rm loop}_i, \bE^{\rm inc})_{L^2_T(\Gamma)}\\
 &= -({\rm curl}_S\, p_i, \bE^{\rm inc})_{L^2_T(\Gamma)}\\
 &= (p_i, {\rm curl}\, \bE^{\rm inc})_{L^2_T(\Gamma)}\\
 &= (p_i, \zi\omega\mu \bH^{\rm inc})_{L^2_T(\Gamma)}
\end{align*}
where $p_i$ is the piecewise linear function introduced in (\ref{eq:def_p}).
If we assume that the incident wave satisfies $\bE^{\rm inc}\sim O(1)$ and $\bH^{\rm inc}\sim O(1)$,
we obtain
\begin{align*}
\bb^L_{L^2} \sim O(k h^3).
\end{align*}
Note that $\phi_i\sim O(h)$ since $\bt_i^{\rm loop}$ is normalised, namely, $O(\bt_i^{\rm loop})\sim 1$.
Also, we obtain
\begin{align*}
\bb^S_{L^2}\sim O(h^2).
\end{align*}
Consequently, the orders of the elements in (\ref{eq:efie_l2_decomposed}) are given as 
\begin{align*}
 \begin{pmatrix}
  Z^{LL}_{L^2}(O(k h^3)) & Z^{LS}_{L^2}(O(k h^3))\\
  Z^{SL}_{L^2}(O(k h^3)) & Z^{SS}_{L^2}\left(O\left(\frac{h}{k} \right)\right)
 \end{pmatrix}
\begin{pmatrix}
 \bj^{\rm loop} \\ \bj^{\rm star}
\end{pmatrix}=
 \begin{pmatrix}
  \bb^L_{L^2}(O(k h^3)) \\ \bb^S_{L^2}(O(h^2))
 \end{pmatrix}.
\end{align*}
Dividing both sides of this equation by $h^2$, we obtain
\begin{align}\label{eq:efie_l2_decomposed_order}
 \begin{pmatrix}
  Z^{LL}_{L^2}(O(k h)) & Z^{LS}_{L^2}(O(k h))\\
  Z^{SL}_{L^2}(O(k h)) & Z^{SS}_{L^2}\left(O\left(\frac{1}{k h} \right)\right)
 \end{pmatrix}
\begin{pmatrix}
 \bj^{\rm loop} \\ \bj^{\rm star}
\end{pmatrix}=
 \begin{pmatrix}
  \bb^L_{L^2}(O(k h)) \\ \bb^S_{L^2}(O(1))
 \end{pmatrix}.
\end{align}
Thus, the orders of all the elements can be written in terms of the powers of $k h$.

We denote the solution of iteration methods after $n$th iteration by
\begin{align}\label{eq:sol_pert}
\begin{pmatrix}
 \bj^{\rm loop}_n \\ \bj^{\rm star}_n
\end{pmatrix}
=
\begin{pmatrix}
 \overline{\bj^{\rm loop}} \\ \overline{\bj^{\rm star}}
\end{pmatrix}
+
\begin{pmatrix}
 \Delta\bj^{\rm loop}_n \\ \Delta\bj^{\rm star}_n
\end{pmatrix}
\end{align}
where $\overline{\bj^{\rm loop}}$ and $\overline{\bj^{\rm star}}$ are the exact solutions of (\ref{eq:efie_l2_decomposed}) and $\Delta\bj^{\rm loop}_n$ and $\Delta\bj^{\rm star}_n$ are the errors of the numerical solutions $\bj^{\rm loop}_n$ and $\bj^{\rm star}_n$.
The orders of the exact solutions $\overline{\bj^{\rm loop}}$ and $\overline{\bj^{\rm star}}$ also can be evaluated from (\ref{eq:efie_l2_decomposed_order}) as
\begin{align}\label{eq:eval_exact_sol}
\overline{\bj^{\rm loop}} \sim O(1),\quad
\overline{\bj^{\rm star}} \sim O(k h).
\end{align}

If the solution satisfies
\begin{align*}
 \frac{\|\bb-A\bx\|}{\|\bb\|} < \delta,
\end{align*}
the iteration method stops and gives (\ref{eq:sol_pert}) as the numerical solution where $\delta$ is the error tolerance.
Substituting (\ref{eq:efie_l2_decomposed_order}) and (\ref{eq:sol_pert}) into this equation, we have
\begin{align*}
\frac{
\left\|
 \begin{pmatrix}
  Z^{LL}_{L^2}(O(k h)) & Z^{LS}_{L^2}(O(k h))\\
  Z^{SL}_{L^2}(O(k h)) & Z^{SS}_{L^2}\left(O\left(\frac{1}{k h}\right)\right)
 \end{pmatrix}
\begin{pmatrix}
 \Delta\bj^{\rm loop}_n \\  \Delta\bj^{\rm star}_n 
\end{pmatrix}
\right\|}
{\left\|
\begin{pmatrix}
\bb_{L^2}^L(O(k h)) \\ \bb_{L^2}^S(O(1))
\end{pmatrix}
\right\|} 
<\delta.
\end{align*}
From this equation and (\ref{eq:eval_exact_sol}), we obtain an estimate of the relative error of the numerical solutions as follows:
\begin{align}\label{eq:l2_res}
 \frac{\|\Delta\bj^{\rm loop}_n\|}{\|\overline{\bj}\|} \sim \frac{\delta}{k h},\quad
 \frac{\|\Delta\bj^{\rm star}_n\|}{\|\overline{\bj}\|} \sim k h \delta.
\end{align}
We thus see that $\bj^{\rm loop}_n$ may have a large relative error if $k h$ is small.


\section{Galerkin Method with the $H_{\rm div}$ Inner Product}\label{sec:disc_proposed}
In this section, we propose a discretisation method,
which achieves good accuracy even in problems with low frequencies.

\subsection{Discretisation}
We utilise the $H_{\rm div}$ inner product
\begin{align*}
 (\bu, \bv)_{H_{\rm div}(\Gamma)} 
:= (\bu, \bv)_{L^2_T(\Gamma)} 
   + c(\nabla_S\cdot \bu, \nabla_S\cdot \bv)_{L^2(\Gamma)}
\end{align*}
for discretising the EFIE in (\ref{eq:efie}) where $c$ is a positive constant.
The constant $c$ is usually set $c=1$ in mathematics.
But we determine the value of $c$ differently in section \ref{sec:hdiv_low_freq} in order to solve the low-frequency breakdown.

For discretisation of (\ref{eq:efie}) with this inner product,
we have to pay attention to the testing function.
First, the testing function should be an $H_{\rm div}$ function
while the $H_{\rm curl}$ function $\bn\times\bt_i$ is used as the testing functions
in the conventional Galerkin method.
Furthermore, we fail if we discretise (\ref{eq:efie}) with the $H_{\rm div}$ inner product
and the RWG testing function as follows:
\begin{align*}
\biggl(\bt_i,\, &\bn\times\int_\Gamma\Bigl\{\zi\omega\mu G(\bx-\by)\bt_j (\by) \\
   &+ \frac{\zi}{\omega\varepsilon}\nabla\nabla G(\bx-\by)\bt_j (\by)\Bigr\}\md S_y\biggr)_{H_{\rm div}(\Gamma)},
\end{align*}
which is ill-conditioned due to the same reason
as the Gram matrix
\begin{align*}
(\bn\times\bt_i,\bt_j)
\end{align*}
is ill-conditioned \cite{christiansen2003preconditioner}.
We can resolve this problem by utilising as a testing function the BC basis function $\bss_i$,
which is the dual function of the RWG function \cite{buffa2005dual}.
We, therefore, construct a discretisation of EFIE with the $H_{\rm div}$ inner product
as follows:
\begin{align}\label{eq:efie_hdiv}
 A_{H_{\rm div}}\bx=\bb_{H_{\rm div}}
\end{align}
where
\begin{align*}
 (A_{H_{\rm div}})_{ij} = \biggl(\bss_i, &\bn\times\int_\Gamma\Bigl\{\zi\omega\mu G(\bx-\by)\bt_j (\by) \\
   &+ \frac{\zi}{\omega\varepsilon}\nabla\nabla G(\bx-\by)\bt_j (\by)\Bigr\}\md S_y\biggr)_{H_{\rm div}(\Gamma)},\\
 (\bb_{H_{\rm div}})_{i} = (\bss_i, &\bE^{\rm inc}\times\bn)_{H_{\rm div}(\Gamma)}.
\end{align*}

The RHS of (\ref{eq:efie_hdiv}) can be calculated as follows:
\begin{align*}
 (\bb_{H_{\rm div}})_{i} &= (\bss_i, \bE^{\rm inc}\times\bn)_{L^2_T(\Gamma)} + 
 c (\nabla_S\cdot\bss_i, \nabla_S\cdot(\bE^{\rm inc}\times\bn))_{L^2(\Gamma)} \\
&= (\bss_i, \bE^{\rm inc}\times\bn)_{L^2_T(\Gamma)} + 
 c (\nabla_S\cdot\bss_i, \bn\cdot(\nabla\times\bE^{\rm inc})_{L^2(\Gamma)} \\
&=(\bss_i, \bE^{\rm inc}\times\bn)_{L^2_T(\Gamma)} + 
 \zi\omega\mu c (\nabla_S\cdot\bss_i, \bn\cdot\bH^{\rm inc})_{L^2(\Gamma)}.
\end{align*}
In a similar way, we can calculate the coefficient matrix whose
second term coincides with the normal component of the magnetic field integral equation (MFIE) as follows:
\begin{align}\notag
 (A_{H_{\rm div}})_{ij} 
 = \biggl(\bss_i, &\bn\times\int_\Gamma\Bigl\{\zi\omega\mu G(\bx-\by)\bt_j (\by) \\\notag
   &+ \frac{\zi}{\omega\varepsilon}\nabla\nabla G(\bx-\by)\bt_j 
(\by)\Bigr\}\md S_y\biggr)_{L^2_T(\Gamma)}\\\label{eq:calc_hdiv_coef}
   - \zi\omega\mu c\Bigl(\nabla_S\cdot\bss_i, &\bn\cdot\int_\Gamma \nabla_y G(\bx-\by)\times\bt_j(\by) 
    \md S_y\Bigr)_{L^2(\Gamma)}.
\end{align}
Hence the discretised integral equation obtained with the $H_{\rm div}$ inner product can be calculated as the sum of the EFIE discretised with the $L^2_T$ inner product and the dual testing functions $\bss_i$, and the normal component of the MFIE tested with the surface divergence of $\bss_i$.
\subsection{Low-Frequency Breakdown}\label{sec:hdiv_low_freq}
We now show that the solution obtained with this discretisation method keeps good accuracy
even in problems with small frequencies.

We apply the loop-star decomposition to the basis functions $\bt_i$ and $\bss_i$ as has been done in section \ref{sec:low_freq_conv}.
Indeed, a linear combination of the BC basis functions can be expanded with the loop and star basis functions  \cite{andriulli2008multiplicative,stephanson2009preconditioned}, namely,
\begin{align*}
 \sum_{i=1}^N c_i \bss_i = \sum_{i=1}^{N_{\rm loop}^{\rm BC}}c_i^{\rm loop} \bss_i^{\rm loop}
 + \sum_{i=1}^{N_{\rm star}^{\rm BC}}c_i^{\rm star} \bss_i^{\rm star}.
\end{align*}
where $N_{\rm loop}^{\rm BC}$ and $N_{\rm star}^{\rm BC}$ are the number of independent loop and star functions of the BC basis functions, respectively.
Note that the loop function $\bss_i^{\rm loop}$ satisfies
\begin{align*}
\nabla_S\cdot \bss_i^{\rm loop} = 0.
\end{align*}
With the help of the loop-star decomposition,  (\ref{eq:efie_hdiv}) reduces to
\begin{align}\label{eq:efie_hdiv_decomposed}
 \begin{pmatrix}
  Z^{LL}_{H_{\rm div}} & Z^{LS}_{H_{\rm div}}\\
  Z^{SL}_{H_{\rm div}} & Z^{SS}_{H_{\rm div}}
 \end{pmatrix}
\begin{pmatrix}
 \bj^{\rm loop} \\ \bj^{\rm star}
\end{pmatrix}=
 \begin{pmatrix}
  \bb^L_{H_{\rm div}} \\ \bb^S_{H_{\rm div}}
 \end{pmatrix}
\end{align}
where $Z^{LL}_{H_{\rm div}}, Z^{LS}_{H_{\rm div}}, Z^{SL}_{H_{\rm div}}$ and $Z^{SS}_{H_{\rm div}}$
 are matrices defined by
\begin{align*}
 (Z^{LL}_{H_{\rm div}})_{ij} &= (\bss^{\rm loop}_i, \zi\omega\mu Q\bt^{\rm loop}_j)_{H_{\rm div}(\Gamma)},\\
 (Z^{LS}_{H_{\rm div}})_{ij} &= (\bss^{\rm loop}_i, \zi\omega\mu Q\bt^{\rm star}_j)_{H_{\rm div}(\Gamma)},\\
 (Z^{SL}_{H_{\rm div}})_{ij} &= (\bss^{\rm star}_i, \zi\omega\mu Q\bt^{\rm loop}_j)_{H_{\rm div}(\Gamma)},\\
  (Z^{SS}_{H_{\rm div}})_{ij} &= (\bss^{\rm star}_i, \zi\omega\mu Q\bt^{\rm star}_j)_{H_{\rm div}(\Gamma)},
\end{align*}
and $\bb^L_{H_{\rm div}}$ and $\bb^S_{H_{\rm div}}$ are vectors defined by:
\begin{align*}
 (\bb^L_{H_{\rm div}})_i &= (\bss^{\rm loop}_i, \bE^{\rm inc}\times\bn)_{H_{\rm div}(\Gamma)},\\
 (\bb^S_{H_{\rm div}})_i &= (\bss^{\rm star}_i, \bE^{\rm inc}\times\bn)_{H_{\rm div}(\Gamma)}.
\end{align*}

We calculate the orders of these elements as in section \ref{sec:low_freq_conv}.
The second term in (\ref{eq:calc_hdiv_coef}) vanishes in $Z^{LL}_{H_{\rm div}}, Z^{LS}_{H_{\rm div}}$ since $\nabla_S\cdot\bss_i^{\rm loop}=0$.
Hence the orders of $Z^{LL}_{H_{\rm div}}, Z^{LS}_{H_{\rm div}}$ are the same
as those of $Z^{SL}_{L^2}, Z^{SS}_{L^2}$. Namely, we have
\begin{align*}
 Z^{LL}_{H_{\rm div}}\sim O(k h^3),\quad Z^{LS}_{H_{\rm div}}\sim O\left(\frac{h}{k}\right).
\end{align*}
Note that $Z^{LS}_{H_{\rm div}}$ has the same order as $Z^{SS}_{L^2}$ since the testing functions of the proposed method do not contain the term $\bn\times$
and, thus, the hyper-singular term in $Z^{LS}_{H_{\rm div}}$ does not vanish.
In $Z^{SL}_{H_{\rm div}}$ and $Z^{SS}_{H_{\rm div}}$,
the second terms do not vanish and their orders depends on the value of the constant $c$.
In $Z^{SL}_{H_{\rm div}}$, for example, we have
\begin{align*}
& (Z^{SL}_{H_{\rm div}})_{ij} \\
&= (\bss^{\rm star}_i, \zi\omega\mu Q\bt^{\rm loop}_j)_{H_{\rm div}(\Gamma)}\\
 &= (\bss^{\rm star}_i, \zi\omega\mu Q\bt^{\rm loop}_j)_{L^2_T(\Gamma)}\\
 &-\zi\omega\mu c \left(\nabla_S\cdot\bss^{\rm star}_i,  \bn\cdot\int_\Gamma \nabla_y G(\bx -\by)\times\bt^{\rm loop}_j\md S_y\right)_{L^2(\Gamma)}\\
  &=\left(\bss^{\rm star}_i, \zi\omega\mu \int_\Gamma G(\bx -\by)\bt^{\rm loop}_j\md S_y\right)_{L^2_T(\Gamma)}\\
 &-\zi\omega\mu c \left(\nabla_S\cdot\bss^{\rm star}_i,  \bn\cdot\int_\Gamma \nabla_y G(\bx -\by)\times\bt^{\rm loop}_j\md S_y\right)_{L^2(\Gamma)}
\end{align*}
since $\nabla_S\cdot \bt^{\rm loop}_j=0$. These two terms satisfy
\begin{align}\label{eq:zsl1}
& \left(\bss^{\rm star}_i, \zi\omega\mu \int_\Gamma G(\bx -\by)\bt^{\rm loop}_j\md S_y\right)_{L^2_T(\Gamma)}\sim O(k h^3), \\\notag
& \zi\omega\mu c \left(\nabla_S\cdot\bss^{\rm star}_i,  \bn\cdot\int_\Gamma \nabla_y G(\bx -\by)\times\bt^{\rm loop}_j\md S_y\right)_{L^2(\Gamma)}\\\label{eq:zsl2}
&\sim O(ck h).
\end{align}
Now, we restrict the value of the constant $c$ in a way that (\ref{eq:zsl2}) is larger than (\ref{eq:zsl1}).
This restriction is identical with the condition
\begin{align}\label{eq:c_cond_temp}
 c>h^2,
\end{align}
 under which we obtain
\begin{align*}
 Z^{SL}_{H_{\rm div}} \sim  O(ck h).
\end{align*}
Similar calculation for $Z^{SS}_{H_{\rm div}}$ yields
\begin{align*}
 Z^{SS}_{H_{\rm div}} \sim  O(ck h)
\end{align*}
under the condition in (\ref{eq:c_cond_temp}).

We can also calculate the RHS as follows:
\begin{align*}
\bb^L_{H_{\rm div}} &= (\bss_i^{\rm loop}, \bE^{\rm inc}\times \bn)_{H_{\rm div}(\Gamma)}\\
 &= (\bss_i^{\rm loop}, \bE^{\rm inc}\times \bn)_{L^2_T(\Gamma)}\\
&\sim O(h^2),\\
 \bb^S_{H_{\rm div}} &=(\bss_i^{\rm star}, \bE^{\rm inc}\times \bn)_{H_{\rm div}(\Gamma)}\\ 
 &= (\bss_i^{\rm star}, \bE^{\rm inc}\times \bn) + \zi\omega\mu c(\nabla_S\cdot \bss_i^{\rm star}, \bn\cdot\bH^{\rm inc})_{L^2(\Gamma)}
\end{align*}
The two terms in $\bb_{H_{\rm div}}^S$ are estimated as
\begin{align}\label{eq:bb_temp1}
&(\bss_i^{\rm star}, \bE^{\rm inc}\times \bn) \sim O(h^2)\\\label{eq:bb_temp2}
&\zi\omega\mu c(\nabla_S\cdot \bss_i^{\rm star}, \bn\cdot\bH^{\rm inc})_{L^2(\Gamma)} \sim O(ck h).
\end{align}
Again, we assume that the term in (\ref{eq:bb_temp2}) including $c$ is larger than the term in (\ref{eq:bb_temp1}),
which leads to 
\begin{align}\label{eq:c_cond_temp2}
c > \frac{h}{k}.
\end{align}
The vector $\bb_{H_{\rm div}}^S$ then satisfies
\begin{align*}
  \bb_{H_{\rm div}}^S \sim O(ck h)
\end{align*}
under the condition in (\ref{eq:c_cond_temp2}).

Hence the orders of elements in (\ref{eq:efie_hdiv_decomposed}) are
 \begin{align*}
 \begin{pmatrix}
  Z^{LL}_{H_{\rm div}}(k h^3) & Z^{LS}_{H_{\rm div}}\left(\frac{h}{k}\right)\\
  Z^{SL}_{H_{\rm div}}(ck h) & Z^{SS}_{H_{\rm div}}(ck h)
 \end{pmatrix}
\begin{pmatrix}
 \bj^{\rm loop} \\ \bj^{\rm star}
\end{pmatrix}=
 \begin{pmatrix}
  \bb^L_{H_{\rm div}}(h^2) \\ \bb^S_{H_{\rm div}}(ck h)
 \end{pmatrix}.
 \end{align*}
 Dividing this equation by $h^2$, we obtain
 \begin{align}\label{eq:efie_hdiv_decomposed_order_tmp}
 \begin{pmatrix}
  Z^{LL}_{H_{\rm div}}(k h) & Z^{LS}_{H_{\rm div}}\left(\frac{1}{k h}\right)\\
  Z^{SL}_{H_{\rm div}}\left(c\frac{k}{h}\right) & Z^{SS}_{H_{\rm div}}\left(c\frac{k}{h}\right)
 \end{pmatrix}
\begin{pmatrix}
 \bj^{\rm loop} \\ \bj^{\rm star}
\end{pmatrix}=
 \begin{pmatrix}
  \bb^L_{H_{\rm div}}(1) \\ \bb^S_{H_{\rm div}}\left(c\frac{k}{h}\right)
 \end{pmatrix}.
 \end{align}
 From this equation, we determine the value of $c$.
 Focusing on the order of $c$ with respect to $k$,
 it is found necessary to take $c=O(1/k^2)$.
 This is because the second row of the coefficient matrix is much larger than the first row
 if $c$ is larger than $O(1/k^2)$,
 and the coefficient matrix approaches a singular matrix when $k\rightarrow 0$.
 If $c$ is smaller than $O(1/k^2)$, the second row in (\ref{eq:efie_hdiv_decomposed_order_tmp}) is much smaller than the first row
 when $k\rightarrow 0$.
 Also $c=O(1/k^2)$ satisfies the restriction in (\ref{eq:c_cond_temp}) and (\ref{eq:c_cond_temp2}) since $|k h|\ll 1$.
 Substituting $c=O(1/k^2)$ in (\ref{eq:efie_hdiv_decomposed_order_tmp}), we have
 \begin{align}\label{eq:efie_hdiv_decomposed_order}
 \begin{pmatrix}
  Z^{LL}_{H_{\rm div}}(k h) & Z^{LS}_{H_{\rm div}}\left(\frac{1}{k h}\right)\\
  Z^{SL}_{H_{\rm div}}\left(\frac{1}{k h}\right) & Z^{SS}_{H_{\rm div}}\left(\frac{1}{k h}\right)
 \end{pmatrix}
\begin{pmatrix}
 \bj^{\rm loop} \\ \bj^{\rm star}
\end{pmatrix}=
 \begin{pmatrix}
  \bb^L_{H_{\rm div}}(1) \\ \bb^S_{H_{\rm div}}\left(\frac{1}{k h}\right)
 \end{pmatrix}.
 \end{align}
 Also from (\ref{eq:efie_hdiv_decomposed_order}), the choice of $c=O(1/k^2)$ seems natural
 since the orders of all the elements are written in terms of powers of $k h$.

 As has been done in section \ref{sec:low_freq_conv},
we decompose the solution $\bj^{\rm loop}$ and $\bj^{\rm star}$
into the exact solution $\overline{\bj^{\rm loop}}$ and $\overline{\bj^{\rm star}}$
and the error $\Delta\bj^{\rm loop}$ and $\Delta\bj^{\rm star}$ respectively.
From (\ref{eq:efie_hdiv_decomposed_order}), we obtain \eqref{eq:eval_exact_sol} again.
The errors $\Delta \bj_n^{\rm loop}, \Delta \bj_n^{\rm star}$ satisfy
\begin{align*}
\frac{
\left\|
 \begin{pmatrix}
  Z^{LL}_{H_{\rm div}}(O(k h)) & Z^{LS}_{H_{\rm div}}\left(O\left(\frac{1}{k h}\right)\right)\\
  Z^{SL}_{H_{\rm div}}(O(\frac{1}{k h})) & Z^{SS}_{H_{\rm div}}\left(O(\frac{1}{k h})\right)
 \end{pmatrix}
\begin{pmatrix}
 \Delta\bj^{\rm loop}_n \\  \Delta\bj^{\rm star}_n 
\end{pmatrix}
\right\|}
{\left\|
 \begin{pmatrix}
  \bb^L_{H_{\rm div}}(1) \\ \bb^S_{H_{\rm div}}\left(\frac{1}{k h}\right)
 \end{pmatrix}
\right\|}<\delta
\end{align*}
for the error tolerance of $\delta$.
This inequality gives error estimates given as follows:
\begin{align*}
 \frac{\|\Delta \bj_n^{\rm loop}\|}{\|\overline{\bj}\|} \sim \delta,\quad 
 \frac{\|\Delta \bj_n^{\rm star}\|}{\|\overline{\bj}\|} \sim \delta.
\end{align*}
We thus conclude that the relative error with this discretisation method is small even if $k h$ is small.


\section{Preconditioning}\label{sec:calderon}
In this section, we discuss preconditioning for the proposed discretisation of EFIE in (\ref{eq:efie_hdiv}). 
In section \ref{sec:calderon_naive}, we introduce a simple Calderon preconditioning which turns out not to be very effective in terms of the computational time.
In section \ref{sec:calderon_helmholtz}, we propose another preconditioning which can successfully reduce the computational time.
 \subsection{Calderon's preconditioning using the single layer potential of Maxwell's equations}\label{sec:calderon_naive}
 In this section, we first introduce a simple way of applying Calderon's preconditioning to the proposed method.
We obtain this preconditioning method by extending the multiplicative Calderon preconditioning \cite{andriulli2008multiplicative} for the conventional EFIE in (\ref{eq:efie_l2}) to the proposed method.
 This preconditioning method is indeed able to decrease the iteration number but is not effective in decreasing the computational time, as we shall see.
 
From Calderon's formulae for Maxwell's equations \cite{nedelec2001acoustic},
we see that the operator $Q$ satisfies
 \begin{align}\label{eq:calderon}
  k^2Q^2 = \frac{\cal I}{4} + {\cal K}
 \end{align}
where ${\cal I}$ is the identity operator and ${\cal K}$ is a compact operator.
This equation implies that the matrix obtained by discretising the operator $Q^2$ is expected to be well-conditioned.
In the conventional Galerkin method, which utilises the $L^2$ inner product as has been shown in section \ref{sec:disc_conv},
$Q^2$ can be discretised into 
\begin{align*}
 A_{L^2}T_{L^2}^{-1}A'_{L^2}T_{L^2}'^{ -1}
\end{align*}
which is known to be well-conditioned \cite{andriulli2008multiplicative} where
\begin{align*}
 (A'_{L^2})_{ij} &= \biggl(\bn\times\bss_i, \bn\times\int_\Gamma\Bigl\{\zi\mu\omega G(\bx-\by)\bss_j (\by) \\
 &+ \frac{\zi}{\omega\varepsilon}\nabla\nabla G(\bx-\by)\bss_j (\by)\Bigr\}\md S_y\biggr)_{L^2_T(\Gamma)},\\
 (T_{L^2})_{ij} &= (\bn\times\bss_i, \bt_j)_{L^2_T(\Gamma)},\quad
 (T'_{L^2})_{ij} = (\bn\times\bt_i, \bss_j)_{L^2_T(\Gamma)}.
\end{align*}
Hence the right preconditioner given by
\begin{align}\label{eq:l2_preconditioner}
 T_{L^2}A'^{-1}_{L-2}T_{L^2}'
\end{align}
is used for solving (\ref{eq:efie_l2}).

This Calderon preconditioning method may appear to be applicable to the $H_{\rm div}$ discretisation since the difference between the conventional and proposed methods is found only in the inner product and the testing function used for the discretisation.
Actually, by discretising (\ref{eq:calderon}) with the $H_{\rm div}$ inner product, we find that the matrix
\begin{align*}
 A_{H_{\rm div}}T^{-1}_{H_{\rm div}}A'_{H_{\rm div}}T'^{-1}_{H_{\rm div}}
\end{align*}
is expected to be well-conditioned where
\begin{align}\notag
 (A'_{H_{\rm div}})_{ij} &= \biggl(\bt_i, \bn\times\int_\Gamma\Bigl\{\zi\mu\omega G(\bx-\by)\bss_j (\by) \\\notag
 &+ \frac{\zi}{\omega\varepsilon}\nabla\nabla G(\bx-\by)\bss_j (\by)\Bigr\}\,\md S_y\biggr)_{H_{\rm div}(\Gamma)},\\\label{eq:hdiv_gram}
 (T_{H_{\rm div}})_{ij} &= (\bt_i, \bt_j)_{H_{\rm div}(\Gamma)},\quad
 (T'_{H_{\rm div}})_{ij} = (\bss_i, \bss_j)_{H_{\rm div}(\Gamma)},
\end{align}
and $A_{H_{\rm div}}$ is the matrix defined in (\ref{eq:calc_hdiv_coef}).
Hence it may seem natural to solve 
\begin{align}\label{eq:hdiv_eq_in_calderon}
 A_{H_{\rm div}} \bx = \bb_{H_{\rm div}}
\end{align}
with the right preconditioner given by
\begin{align}\label{eq:hdiv_naive_preconditioner}
 T_{H_{\rm div}}A_{H_{\rm div}}'^{-1}T_{H_{\rm div}}'.
\end{align}
This preconditioning indeed decreases the iteration number of linear solvers for
(\ref{eq:hdiv_eq_in_calderon})
but the whole computational time may not be reduced efficiently
since the Gram matrices $T_{H_{\rm div}}$ and $T_{H_{\rm div}}'$ are ill-conditioned.
In fact, the Gram matrix $T_{H_{\rm div}}$ is written as
\begin{align*}
 (T_{H_{\rm div}})_{ij} &= (\bt_i, \bt_j)_{L^2_T(\Gamma)}
 + c(\nabla_S\cdot\bt_i, \nabla_S\cdot \bt_j)_{L^2(\Gamma)}
\end{align*}
with
\begin{align*}
 (\bt_i, \bt_j)_{L^2_T(\Gamma)} \sim O(h^2),\quad c(\nabla_S\cdot\bt_i, \nabla_S\cdot \bt_j)_{L^2(\Gamma)} \sim O(c).
\end{align*}
Thus, the choice $c=O(1/k^2)$ obtained in section \ref{sec:hdiv_low_freq} 
gives
\begin{align*}
 (T_{H_{\rm div}})_{ij} \sim \frac{1}{k^2}(\nabla_S\cdot\bt_i, \nabla_S\cdot \bt_j)_{L^2(\Gamma)}
\end{align*}
as $k h\rightarrow 0$,
which is a singular matrix.
The matrix $T_{H_{\rm div}}'$ also has the same ill-conditioning.
Hence the Gram matrices $T_{H_{\rm div}}$ and $T_{H_{\rm div}}'$ are ill-conditioned in the low frequency region.
This causes much computational time to invert these Gram matrices and, even worse, the failure of the preconditioning for smaller frequencies as will be shown in section \ref{sec:num_ex}.

\subsection{Preconditioning Using Single Layer Potential of Helmholtz' Equation}\label{sec:calderon_helmholtz}
We propose a new preconditioning for the  $H_{\rm div}$-inner-product-discretised EFIE in this section.
This preconditioning will be shown to decrease the iteration number efficiently
and the related Gram matrices to be well-conditioned.

We first note that the coefficient matrix in (\ref{eq:calc_hdiv_coef}) can be written as in (\ref{eq:calc_hdiv_mat}).
\begin{table*}[!htbp] \centering
\begin{align}\notag
  (A_{H_{\rm div}})_{ij} 
 =& \left(\bss_i, \bn\times\int_\Gamma\left\{\zi\omega\mu G(\bx-\by)\bt_j (\by) 
   + \frac{\zi}{\omega\varepsilon}\nabla\nabla G(\bx-\by)\bt_j 
(\by)\right\}\,\md S_y\right)_{L^2_T(\Gamma)}
   - \zi\omega\mu c\left(\nabla_S\cdot\bss_i, \bn\cdot\int_\Gamma \nabla_y G(\bx-\by)\times\bt_j(\by) 
 \,\md S_y\right)_{L^2(\Gamma)} \\\notag
 =& \left(\bss_i, \bn\times\int_\Gamma\left\{\zi\omega\mu G(\bx-\by)\bt_j (\by) 
   + \frac{\zi}{\omega\varepsilon}\nabla\nabla G(\bx-\by)\bt_j 
(\by)\right\}\,\md S_y\right)_{L^2_T(\Gamma)}
   + \zi\omega\mu c\left(\bss_i, \nabla_S\, \bn\cdot\int_\Gamma \nabla_y G(\bx-\by)\times\bt_j(\by) 
 \,\md S_y\right)_{L^2(\Gamma)} \\\label{eq:calc_hdiv_mat}
 =& \biggl(\bss_i, \bn\times\int_\Gamma\left\{\zi\omega\mu G(\bx-\by)\bt_j (\by) 
 + \frac{\zi}{\omega\varepsilon}\nabla\nabla G(\bx-\by)\bt_j
 (\by)\right\}\,\md S_y
 + \zi\omega\mu c \,\nabla_S\, \bn\cdot\int_\Gamma \nabla_y G(\bx-\by)\times\bt_j(\by) 
 \,\md S_y
 \biggr)_{L^2_T(\Gamma)}
\end{align}
\end{table*}
Hence the coefficient matrix $A_{H_{\rm div}}$ can be regarded as the matrix obtained by discretising the integral operator
\begin{align*}
 &\widetilde{Q}\bu \\
=& \bn\times\int_\Gamma\left\{\zi\omega\mu G(\bx-\by)\bu (\by) 
 + \frac{\zi}{\omega\varepsilon}\nabla\nabla G(\bx-\by)\bu
 (\by)\right\}\md S_y\\
 +& \zi\omega\mu c \nabla_S \bn\cdot\int_\Gamma \nabla_y G(\bx-\by)\times\bu(\by) 
 \md S_y
\end{align*}
with the $L^2$ inner product and the testing function $\bss_i$.

Now we construct a preconditioner for the operator $\widetilde{Q}$ with the help of principal symbols.
We take a local coordinate in the tangential plane on the boundary $\Gamma$ whose 3rd axis is directed in the direction of the normal vector $\bn$.
We then compute the Fourier transforms of the singular parts of the integral operator $\widetilde{Q}$ within the tangential plane. The result is
\begin{align}\label{eq:principal_q}
 \frac{\zi\omega\mu}{2}
 \left(-\frac{\epsilon_{ij}}{\rho}
 + \frac{\epsilon_{ip}\xi_p\xi_j}{k^2\rho}
 + c \frac{\xi_i\epsilon_{jp}\xi_p}{\rho}\right)
\end{align}
where $\epsilon_{ij}$ is the permutation symbol in 2D, $\xi_i\, (i=1,2)$ is the Fourier parameter and $\rho=\sqrt{|\xi|^2-k^2}$. Note that we use the summation convention to repeated indices in this equation as well as in the rest of this section.
We next introduce an operator $\widetilde{\cal S}$ defined by
\begin{align*}
\zi\omega\varepsilon \widetilde{\cal S}\bu = \zi\omega\varepsilon\bn\times \int_\Gamma G(\bx - \by)\bu(\by)\, \md S_y,
\end{align*}
which is included as a part in the operator $\widetilde{Q}$.
This operator $\widetilde{\cal S}$ has a principal symbol given by
\begin{align}\label{eq:principal_s}
-\frac{\zi\omega\varepsilon\epsilon_{jk}}{2\rho}.
\end{align}
Hence the product of (\ref{eq:principal_q}) and (\ref{eq:principal_s}) is asymptotically equal to
 \begin{align}\label{eq:principal_qs}
  p_0(\widetilde{Q}\cdot\zi\omega\varepsilon\widetilde{\cal S})=
 -\frac{k^2}{4|\bxi|^2}\left(
 \frac{1}{k^2}\epsilon_{ip}\xi_p\epsilon_{kj}\xi_j
 - c\xi_i\xi_k
 \right)
 \end{align}
 as $|\xi| \rightarrow \infty$.
 The matrix $p_0(\widetilde{Q}\cdot\zi\omega\varepsilon\widetilde{\cal S})$, or the principal symbol of the operator $\widetilde{Q}\cdot\zi\omega\varepsilon\widetilde{\cal S}$, determines the operator $\widetilde{Q}\cdot\zi\omega\varepsilon\widetilde{\cal S}$ to within a compact operator.
 The eigenvectors of this matrix are obviously $\xi_k$ and $\epsilon_{kq}\xi_q$,
 and their eigenvalues are $ck^2/4$ and $-1/4$, respectively.
 Thus we conclude that
 \begin{align*}
  \widetilde{Q}\cdot\zi\omega\varepsilon \widetilde{\cal S} = {\cal R} + {\cal K}
 \end{align*}
 where ${\cal R}$ is an operator on $\Gamma$ whose eigenvalues are $ck^2/4$ and $-1/4$
 and ${\cal K}$ is a compact operator.
 In other words, the eigenvalues of the operator $\widetilde{Q}\cdot\zi\omega\varepsilon \widetilde{\cal S}$ accumulate at $ck^2/4$ and $-1/4$.
 The operator $\widetilde{Q}\cdot\zi\omega\varepsilon \widetilde{\cal S}$ is discretised into
 \begin{align}\label{eq:hdiv_new_calderon_mat}
  A_{H_{\rm div}} T_{L^2}^{-1}\widetilde{S}_{L^2}T''^{-1}_{L^2}
 \end{align}
 where
 \begin{align*}
  (A_{H_{\rm div}})_{ij} &= (\bss_i, \zi\omega\mu Q\bt_j)_{H_{\rm div}(\Gamma)}
  = (\bss_i, \zi\omega\mu \widetilde{Q}\bt_j)_{L^2(\Gamma)},\\
  (\widetilde{S}_{L^2})_{ij} &= (\bn\times\bss_i, \zi\omega\varepsilon\widetilde{S}\bss_j)_{L^2(\Gamma)},\\
  (T''_{L^2})_{ij} &= (\bss_i, \bss_j)_{L^2(\Gamma)}.
 \end{align*}
Note that the operator $\widetilde{\cal Q}$ is introduced only for the explanation of the preconditioning based on $\widetilde{S}_{L^2}$ but is never used in computation.

  Consequently, we find that the eigenvalues of the matrix in (\ref{eq:hdiv_new_calderon_mat}) are expected to accumulate around $ck^2/4$ and $-1/4$. 
In section \ref{sec:hdiv_low_freq}, we found that the low-frequency breakdown can be solved with $c=O(1/k^2)$.
This choice of $c$ is also suitable for this preconditioning
since the condition number of the matrix in (\ref{eq:principal_qs}) is bounded and even becomes $1$ with $c=1/k^2$.
As a result, (\ref{eq:efie_hdiv}) can be preconditioned with the following right preconditioner:
\begin{align}\label{eq:hdiv_new_calderon_preconditioner}
 T_{L^2}\widetilde{S}_{L^2}^{-1}T''_{L^2}
\end{align} 
with $c=1/k^2$.
To use this preconditioner, we need inversion of the matrices $T_{L^2}$ and $T''_{L^2}$.
These inversions, however, do not take much computational time since
the Gram matrices $T_{L^2}$ and $T''_{L^2}$ are well-conditioned in contrast to $T_{H_{\rm div}}$ and $T'_{H_{\rm div}}$, which appear in the preconditioner in (\ref{eq:hdiv_naive_preconditioner}).
We note, however, that the use of $\widetilde{S}_{L^2}$ for the preconditioner may cause spurious resonances  in addition to those of the EFIE, although $\widetilde{S}_{L^2}$ is otherwise a regular matrix.

 \section{Numerical Examples}\label{sec:num_ex}
 The following five different combinations of the discretisation methods and the preconditioning methods are tested in this section.
 \begin{itemize}
  \item Approach 1: The $H_{\rm div}$ inner product with the preconditioning proposed in section \ref{sec:calderon_helmholtz} ((\ref{eq:efie_hdiv}) is solved with the right preconditioner in (\ref{eq:hdiv_new_calderon_preconditioner})).
  \item Approach 2: The $H_{\rm div}$ inner product with the preconditioning proposed in section \ref{sec:calderon_naive} ((\ref{eq:efie_hdiv}) is solved with the right preconditioner in (\ref{eq:hdiv_naive_preconditioner})).
  \item Approach 3: The $L^2$ inner product with the Calderon preconditioning ((\ref{eq:efie_l2}) is solved with the right preconditioner in (\ref{eq:l2_preconditioner})).
  \item Approach 4: The $H_{\rm div}$ inner product without preconditionings ((\ref{eq:efie_hdiv}) is solved without preconditioning).
  \item Approach 5: The $L^2$ inner product without preconditionings ((\ref{eq:efie_l2}) is solved without preconditioning).
 \end{itemize}
In our implementation, we compute hypersingular integrals in the matrices in \eqref{eq:efie_l2} and \eqref{eq:efie_hdiv} after regularisation using integration by parts.
Both derivatives in $\nabla\nabla G$ are moved to trial functions in \eqref{eq:efie_l2}
while only one of the derivatives are moved in \eqref{eq:efie_hdiv}.

\subsection{Spherical Scatterer}
A spherical PEC with the radius of $0.25$ illuminated by the plane incident wave given by
\begin{align*}
 \bE^{\rm inc}(\bx) = \bE^{\rm inc}_0\ex^{\zi\bk\cdot\bx}
\end{align*}
is considered where
\begin{align*}
 \bk=(0,0,k)^T,\quad \bE^{\rm inc}_0 = (1,0,0)^T.
\end{align*}
We set $\varepsilon=\mu=1$ in the exterior domain $\Omega^e$.
The frequency is nondimensionalised such that
the wavelength is equal to one when the frequency $k$ is $2\pi$.
The surface of the spherical scatterer is divided with the meshes with 10580 and 128000 triangular elements.
The RWG and BC basis functions are used for $\bt_i$ and $\bss_i$, respectively.
The GMRES with the error tolerance of $10^{-5}$ is used for both solving the discretised integral equation and calculating the inverse of the Gram matrices.
The low-frequency FMM is used for accelerating the computation of the coefficient matrix.

\begin{figure}[htbp]
  \centering
      \includegraphics[clip, width=8cm]{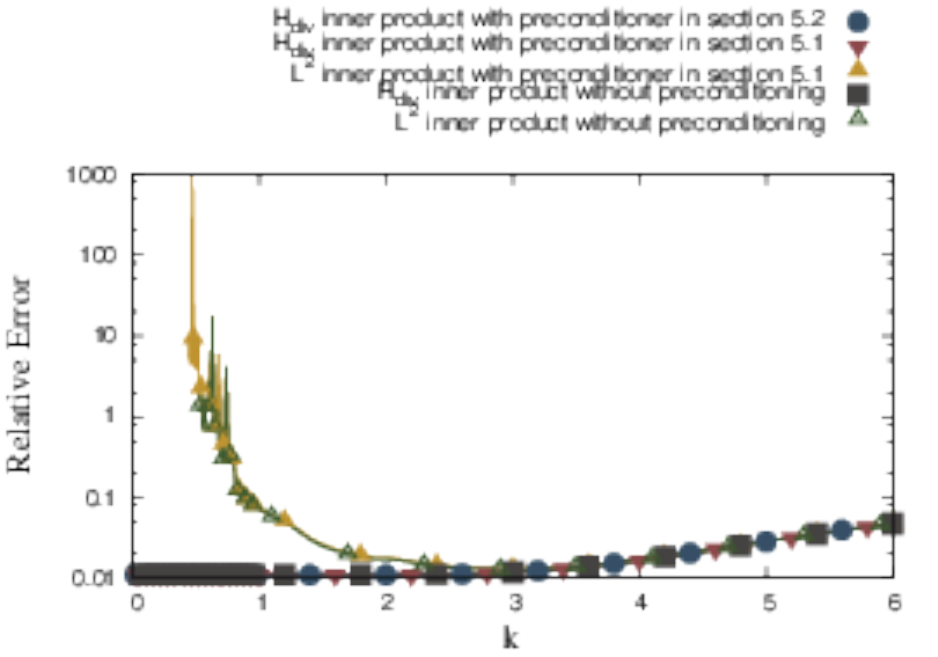}
    \caption{The relative error with 10580 triangular elements.
The yellow and green lines are truncated since the GMRES does not converge after 3000 iterations.}
    \label{omg_coarse_err}
\end{figure}

\begin{figure}[htbp]
  \centering
      \includegraphics[clip, width=8cm]{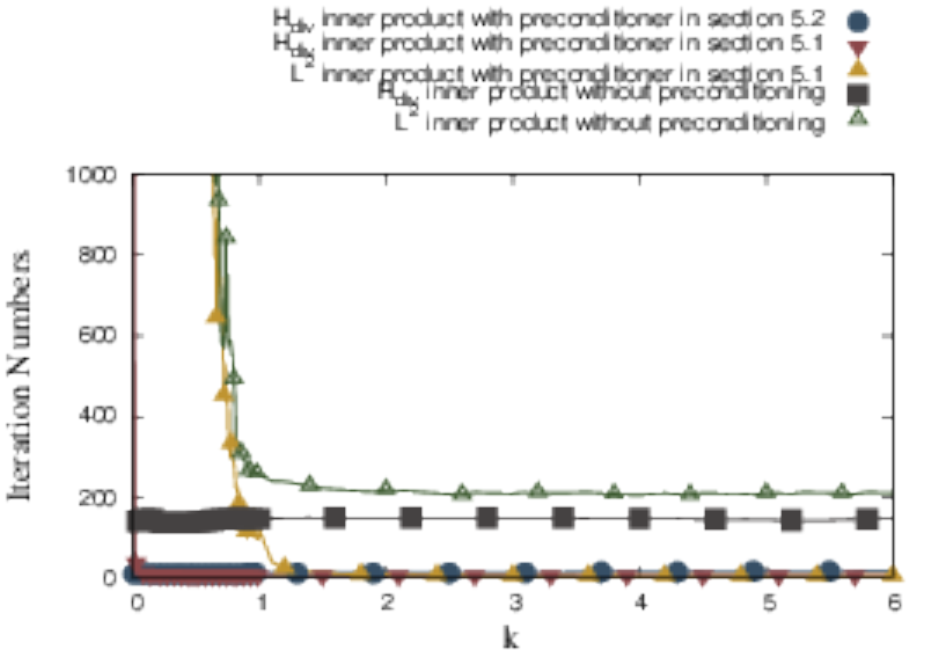}
    \caption{The iteration number of the GMRES with 10580 triangular elements.
}
    \label{omg_coarse_iter}
\end{figure}

\begin{figure}[htbp]
  \centering
      \includegraphics[clip, width=8cm]{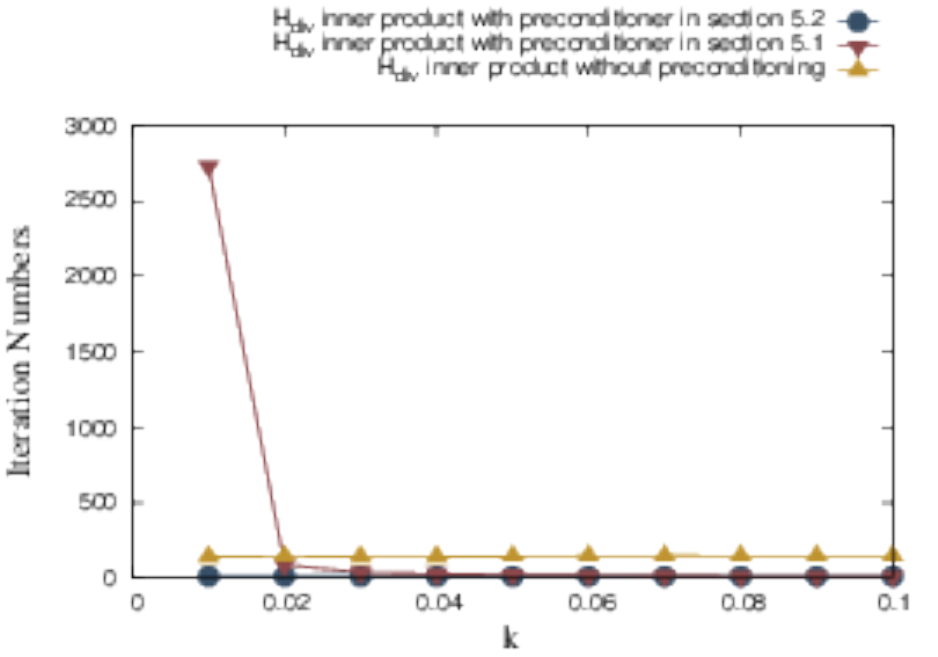}
    \caption{The iteration number of the GMRES with 10580 triangular elements.}
    \label{omg_coarse_iter_small}
\end{figure}

\figuref\ref{omg_coarse_err} shows the relative error
of the numerical methods for the mesh with 10580 triangular elements.
The relative error is defined by
\begin{align*}
 \frac{\sqrt{\int_\Gamma \|\bj_{\rm cal} - \bj_{\rm ana}\|^2\md S}}{\sqrt{\int_\Gamma \|\bj_{\rm ana}\|^2\md S}}
\end{align*}
where $\bj_{\rm cal}$ is the numerical solution and $\bj_{\rm ana}$ is the analytic solution obtained with the Mie series.
The yellow and green lines (approaches 3 and 5) in \figuref\ref{omg_coarse_err} are truncated since we set the maximum iteration number of the GMRES to be $3000$ in this example and the GMRES in approaches 3 and 5 did not converge after the maximum iterations in some small frequencies.
The methods with the $H_{\rm div}$ inner product (approaches 1, 2 and 4) show good accuracy for any frequency
while the accuracy of the methods with the $L^2$ inner product (approaches 3 and 5) becomes worse as the frequency decreases.
The relative errors of the three methods using the $H_{\rm div}$ inner product
are almost the same.
This implies that the relative error is independ of the preconditioning methods, as it should be.
\figuref\ref{omg_coarse_iter} shows the iteration number of the GMRES for the same problem.
The methods with preconditioning (approaches 1 $\sim$ 3) require much less iteration numbers than those without preconditioning (approaches 4 and 5).
The iteration number with the $L^2$ inner product (approaches 3 and 5) diverges when $k < 1$
since the coefficient matrices of these methods are almost singular in these frequencies.
\figuref\ref{omg_coarse_iter_small} also shows the iteration number of the methods using the $H_{\rm div}$ inner product but the region of the frequency $k$ is restricted to $0.01< k < 0.1$.
From \figuref\ref{omg_coarse_iter_small}, we find that the iteration number of approach 2 increases for very small frequencies ($k \sim 0.01$).
This is because the Gram matrices are ill-conditioned for small frequencies as stated in section \ref{sec:calderon_naive}.

\begin{figure}[htbp]
  \centering
      \includegraphics[clip, width=8cm]{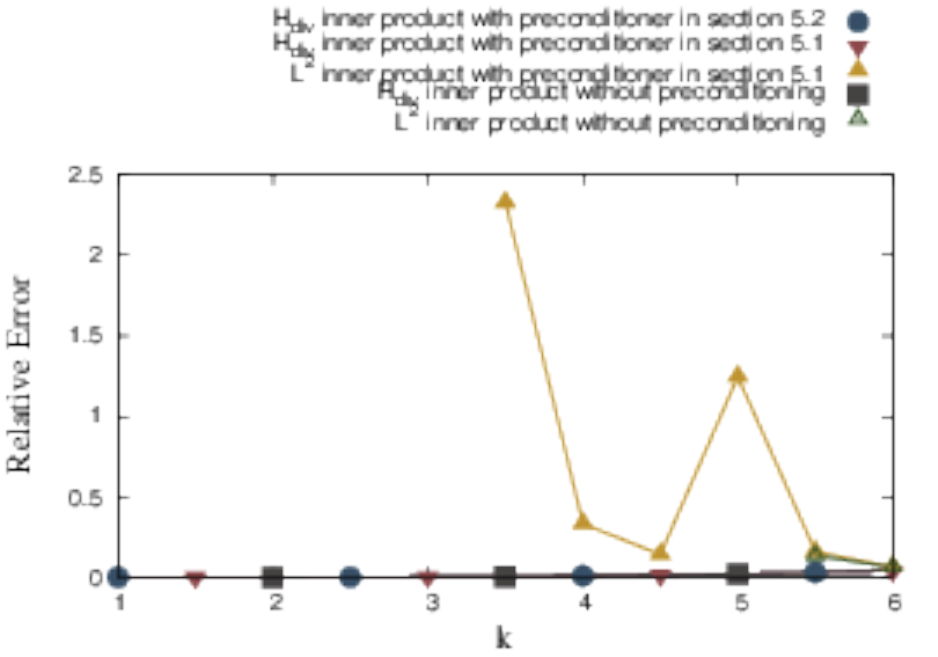}
 \caption{The relative error with 128000 triangular elements.
The yellow and green lines are truncated since the GMRES does not converge after 1500 iterations.}
    \label{omg_fine_err}
\end{figure}

\begin{figure}[htbp]
  \centering
      \includegraphics[clip, width=8cm]{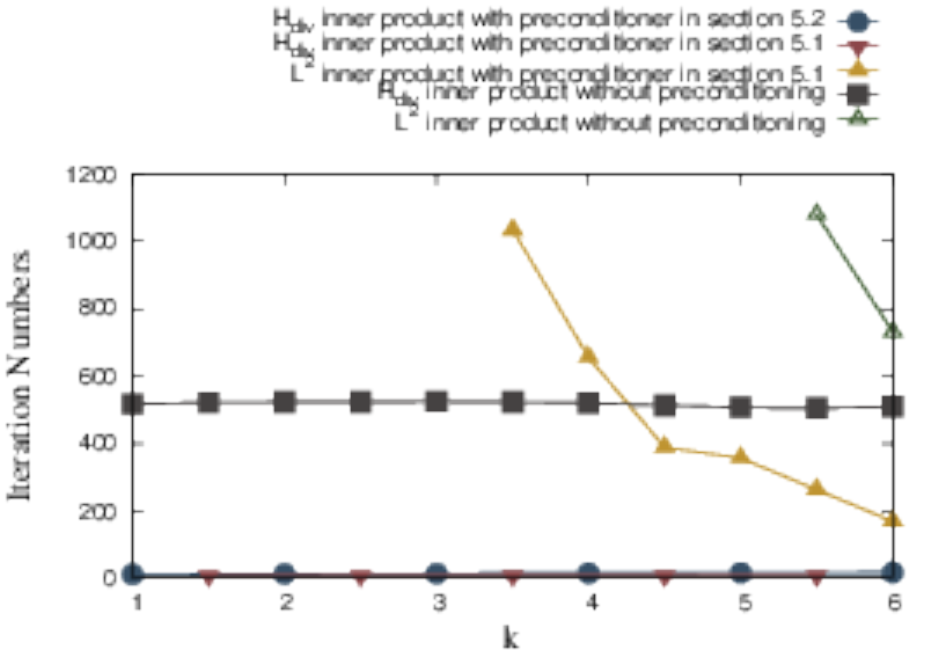}
 \caption{The iteration number of the GMRES with 128000 triangular elements.
 The yellow and green lines are truncated since the GMRES does not converge after 1500 iterations.
}
    \label{omg_fine_iter}
\end{figure}

\begin{figure}[htbp]
  \centering
      \includegraphics[clip, width=8cm]{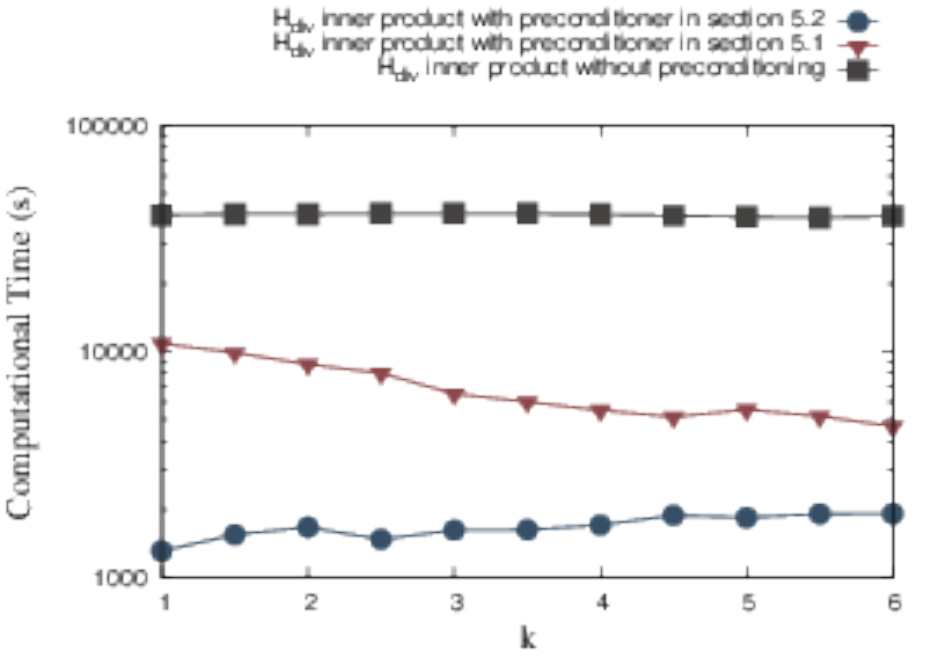}
 \caption{The computational time of the GMRES with 128000 triangular elements.
}
    \label{omg_fine_time}
\end{figure}

\figuref\ref{omg_fine_err}
shows the relative error for the finer mesh with 128000 triangular elements. 
We set the maximum iteration numbers of the GMRES to be $3000$ in this example.
The lines for approaches 3 and 5 in this figure are again truncated
since the GMRES after the maximum iterations did not reach the error tolerance at the omitted points.
In this example, $kh$ is smaller for the same $k$ than that in \figuref\ref{omg_coarse_err} since the mesh size $h$ is smaller.
Hence, in the methods with the $L^2$ inner product (approaches 3 and 5), the relative error is larger than the results in \figuref\ref{omg_coarse_err} or the GMRES did not converge for almost all frequencies in \figuref\ref{omg_fine_err}.
Methods with the $H_{\rm div}$ inner product, however, show good accuracy even for such a fine mesh.
\figuref\ref{omg_fine_iter} shows the iteration number for the same example.
The methods with the $L^2$ inner product required a large number of iterations and
 did not
reach the error tolerance after the maximum iteration number of $3000$ in many cases.
Comparing the three methods based on the $H_{\rm div}$ inner product, we see that the combinations of the $H_{\rm div}$ inner product with the preconditionings (approaches 1 and 2) lead to convergence with about ten iterations while the method without preconditioning (approach 4) requires about 500 iterations. 
\figuref\ref{omg_fine_time} shows the computational time of approaches 1, 2 and 4, which are based on the $H_{\rm div}$ inner product.
The computational time of approach 2 is much more than that of approach 1 and increases as the frequency goes smaller
even though the iteration numbers of approaches 1 and 2 are almost the same.
This is due to the inversion of the ill-conditioned Gram matrices in (\ref{eq:hdiv_gram}) in approach 2, which is stated in section \ref{sec:calderon_naive}.
In fact, as shown in \tableref\ref{omg_fine_avetime}, the average computational time for a matrix vector product is not different in approaches 1 and 2
but the inversion of the Gram matrices in approach 2 requires much more computational time than that in approach 1 when $k=1$. 
From this result, we conclude that approach 1 is better than approach 2 in terms of the computational time.

 \begin{table}[htbp]
  \caption{Average Computational Time (s) for a Matrix-Vector Product and an Inversion of the Gram Matrices When $k=1$}
  \label{omg_fine_avetime}
 \begin{tabular}{|c|c|c|}\hline
  & product of the matrix (\ref{eq:efie_hdiv})&  inversion of the Gram matrices\\\hline
  approach 1 & 77.34 & 28.08\\\hline
  approach 2 & 77.27 & 1253.98\\\hline
 \end{tabular}
 \end{table}

\section{Conclusion}\label{sec:conclusion}
We proposed a Galerkin method with the $H_{\rm div}$ inner product.
This discretisation method resolves the low-frequency breakdown of the EFIE.
We also described two preconditioners for this method,
one based on the Calderon's formula and another using a part of the EFIO. We have verified that the latter preconditioning using the matrix in (\ref{eq:hdiv_new_calderon_preconditioner}) as a right preconditioner is better in terms of the computational time than the Calderon preconditioner, although the Calderon preconditioner could also reduce the iteration number.

In this paper, we have tested the proposed method in simple problems with small frequencies
in order to make sure that it resolves the low-frequency breakdown.
The behaviours of the proposed method in problems with scatterers of complicated shapes or with higher frequencies, however, remain to be investigated.
Also, we did not deal with spurious resonances in this paper, including
those introduced possibly by the preconditioning operator $\widetilde{\cal S}$, which is another remaining issue. But we expect that the latter problem can be resolved with the help of methods of ``complexified'' wave number \cite{contopanagos2002well}
or simply by taking $k=0$ in $\widetilde{S}$.


%



\section*{Acknowledgment}
This work is supported by JSPS KAKENHI Grant Number 26790078.

\ifCLASSOPTIONcaptionsoff
  \newpage
\fi



\bibliographystyle{IEEEtran}
\bibliography{paper_en}
\end{document}